\documentclass[11pt]{article}
\usepackage{amsfonts}
\usepackage{amssymb}
\usepackage{amsmath}
\usepackage{epsfig,epic,eepic}

\usepackage{psfrag}

\newcommand{\R}{\mathbb R}
\newcommand{\N}{\mathbb N}

\newcommand{\C}{\mathbb C}

\newcommand{\hop}{\vskip .2cm\noindent}
\newcommand{\hip}{\vskip .1cm\noindent}

\newcommand{\mrm}{\mathrm}
\newcommand{\hD}{\widehat D}
\newcommand{\ff}{\mathcal F}

\newcommand{\hM}{\widehat M}
\newcommand{\hg}{\widehat g}

\newcommand{\wX}{X}
\newcommand{\wY}{Y}

\renewcommand{\H}{\mathbb H}

\newtheorem{enonce}{}[section]

\newtheorem{thm}{Theorem}
\newtheorem{cor}[enonce]{Corollary}
\newtheorem{prop}[enonce]{Proposition}
\newtheorem{lem}[enonce]{Lemma}
\newtheorem{defi}[enonce]{Definition}
\newtheorem{rema}[enonce]{Remark}
\newtheorem{fact}[enonce]{Fact}

\newtheorem{example}[enonce]{Example}

\author {Pierre Mounoud}
\title{On parallel and symmetric $2$-tensorfields on cones over pseudo-Riemannian manifolds.}
\date{}
\begin{document}
\maketitle
\begin{abstract}
In this article,  we study complete pseudo-Riemannian manifolds whose cone admits a parallel symmetric $2$-tensorfield.  The situation splits in three cases: nilpotent, decomposable or complex Riemannian. In the complex Riemannian and decomposable cases we provide a classification. In the nilpotent case, we are able to describe completely only a dense open subset of the manifold. To conclude, we give examples with non-constant curvature in the nilpotent case.
\end{abstract}
\section{Introduction.}
Let $(M,g)$ be a pseudo-Riemannian manifold (we will consider that a Riemannian metric is also a pseudo-Riemannian one). It is interesting to associate to $(M,g)$ its \emph{cone} $(\hM,\hg)$ defined by $\hM=\R_+^*\times M$ and $\hg=dr^2+r^2g$. 
This construction appears in different contexts. For example (pseudo-)Sasakian manifolds 
are characterized by the fact that their cone is (pseudo-)K\"alher (see \cite{bookSasakian}).  Cones are also used by B\"ar in \cite {Bar} to classify the Killing spinor on Riemannian manifolds: indeed a Riemannian manifold admits a real Killing spinor if and only if its cone admits a parallel spinor (see \cite{Boh} for the pseudo-Riemannian analogue).  
In \cite{Ga}, Gallot used cones 
to study the Obata equation $(*)$. 
This situation led Alekseevski and al. to study in \cite{Leist} the holonomy of cones as a subject on his own.


In this paper we are interested in describing the pseudo-Riemannian manifolds $(M,g)$ whose  cone has non trivial parallel symmetric $2$-tensor. As such it can be seen as a contribution to the study of the holonomy of cones, in the same spirit as \cite{Leist}. Indeed most of our results extend,  more or less directly,  results of \cite{Leist}. Beside this aspect, there is another reason to be interested in cones admitting such tensors.
Gallot proved in \cite{Ga} (see also \cite{Matmoun}) that they are in one-to-one correspondence with functions $\alpha$ on $M$ which are solutions of the following Obata equation:
$$(*)\qquad DDD \alpha (X,Y,Z)+ 2 (D\alpha  \otimes g)(X,Y,Z)+(D\alpha \otimes g)(Y,X,Z)+(D\alpha \otimes g)(Z,X,Y)=0,$$
where $X$, $Y$, $Z$ are vectors tangent $M$ and $D$ is the Levi-Civita connection of $g$. 

The interest of Gallot in this equation was coming from spectral geometry as the third eigenvalue of the Laplacian of the  sphere is a solution of $(*)$. But this equation is present in projective geometry: it appears in the work of Solodovnikov \cite{So} (see also  the work of Matveev and Kiosak \cite{KM} in the pseudo-Riemannian case) where it is related to the existence of projectively equivalent metrics ie metrics having the same unparameterized geodesics.

 Gallot and Tanno have independently shown (see \cite{Ga} and \cite{Ta})  that if $(M,g)$ is a Riemannian complete manifold admitting a non-constant solution to equation $(*)$ then $(M,g)$ is a quotient of the round sphere. This result has been extended to the pseudo-Riemannian case, keeping the conclusion, by Matveev and the author in \cite{Matmoun} under the hypothesis that $M$ is compact but not necessarily complete. 

Hence  this article can also be seen as an attempt to obtain a more general version of  Gallot-Tanno Theorem. 
In order to obtain global results the hypothesis of compactness used in \cite{Matmoun} has to be replaced by something. 
The first idea that comes in mind is to use the original hypothesis of geodesical completeness.  But, in the pseudo-Riemannian context this hypothesis is not as natural as it is in Riemannian geometry. Thus, we will also consider the  alternative hypothesis 
that there exists a solution of $(*)$ which is proper, ie such that for any compact $K\subset \R$, $\alpha^{-1}(K)$ is compact. 
\vskip .3cm
Let us consider a parallel tensorfield $T$ on $\hM$ and $\widetilde T$ the endomorphism associated to $T$ (ie we have $T(u,v)=\hg(u,\widetilde T (v))$). There are actually $3$ cases to consider (see Proposition \ref{above}):
\begin{itemize}
\item when $T$ is $2$ step nilpotent ie when $\widetilde T^2=0$
\item when $(\hM,\hg,T)$ defines complex Riemannian structure ie when $\widetilde T^2=-\mrm{Id}$
\item when $(\hM,\hg)$ is decomposable ie when $\widetilde T^2=\widetilde T$
\end{itemize}
Our first result is a complete description of an open dense set of $M$. It turns out that, in the complex Riemannian case this open set is $M$ itself, hence this case is understood and we already have lots of manifolds admiting non trivial solutions to equation $(*)$.
\hop
{\bf Theorem \ref{lieuregulier}}
{\sl Let $(M,g)$ be a pseudo-Riemannian manifold  such that $(\hM,\hg)$ has a parallel symmetric $2$ tensorfield $T$. 
 If $(M,g)$ is complete or if $\alpha$ has compact levels then  there exist an open dense set $O\in M$, an open interval $I$, a pseudo-Riemannian manifold $(N,h)$ endowed with a parallel symmetric $2$ tensorfield $S$ 
such that if $U$ is a connected component of $O$ then $U=I\times N$ and
\begin{enumerate}
 \item 
If $\widetilde T^2=0$, then 
$$ 
 g=-ds^2 +e^{2s}(h-S)+S.
$$
\item 
$\widetilde T^2=-\mrm{Id}$, then  
$$ 
 g=-ds^2 + h
- \sinh(2s) S,
$$
Moreover $U=M$.
\item 
If $\widetilde T^2=\widetilde T$ 
$$ 
\begin{array}{ll}
& g=\sin^2(s) (h-S) + \cos^2(s) S\\
\mrm{or}  &g=\sinh^2(s) (h-S) + \cosh^2(s) S.
\end{array}
$$
\end{enumerate}
Reciprocally, the cone over any of those manifolds $(I\times N,g)$ admits a parallel symmetric $2$ tensorfield $T$.}
\hop
The third point and the special case $S=0$ of the first point of Theorem \ref{lieuregulier} were already proven in \cite{Leist}, where they appear respectively as a part of Theorem 7.1 and as Theorem 9.1.

The second important result is the classification, under any of our hypothesis of completeness, of pseudo-Riemannian manifolds with decomposable cone. It extends Theorem 7.1 of \cite{Leist}.
\hop
{\bf Theorem \ref{cas2} and \ref{cas3}}
{\sl Let $(M,g)$ be a pseudo-Riemannian manifold with decomposable cone. 
If $g$ is \emph{complete} or $\alpha$ is proper,
then 
\begin{itemize}
\item either $g$ has constant curvature equal to $1$ 
\item or 
$\widetilde M$, the universal cover of $M$, is a warped product of a negative hyperbolic space (possibly $1$ dimensional) and 
a pseudo-Riemannian manifold $(N,h)$. 
\end{itemize}}
Hence, there is only one case left: when the cone has a nilpotent parallel endomorphism. Unfortunately, we have not been able to give a classification in this case. Instead we provide in section \ref{sect-nilpo} a family of examples consisting of perturbation of the pseudo-sphere $S^{p,q}=\{x\in \R^{p+1,q}\,|\, \langle x,x\rangle=1\}$ where $\R^{p+1,q}$ stand for $\R^{p+q+1}$ endowed with a quadratic form of signature $(p+1,q)$. 
It is our opinion that this lack of rigidity explains why the classification is more complicated in this case.

We conclude this paper, by expliciting how in certain condition a parallel $2$-tensor on $\hM$ provide  projectively equivalent metrics on $M$.
\section{Parallel symmetric $2$-tensors on the cone over a manifold.}\label{section}
\subsection{Link with the Obata equation.}
We start by giving the definition of cones over pseudo-Riemannian manifolds.
\begin{defi}
Let $(M,g)$ be a pseudo-Riemannian manifold. We call \emph{cone manifold} over $(M,g)$ the manifold $\hM=\R^*_+\times M$ endowed with the metric $\hg$ defined by $\hg=\,dr^2+ r^2 g$.
\end{defi}
We will denote by $D$ the Levi-Civita connection of $g$ and by $\hD$ the Levi-Civita connection of $\hg$. Those connections are related by the following fact.
\begin{fact}\label{fact}
The Levi-Civita connection of $\hg$ is given by 
$$\hD_{\wX}\wY={D_XY}-rg(X,Y)\partial_r,\quad \hD_{\partial_r}\partial_r=0, \quad \hD_{\partial_r}\wX=\hD_{\wX}\partial_r=\frac{1}{r}\wX.$$
\end{fact}
 The holonomy of cones over pseudo-Riemannian is strongly related to the equation $(*)$ seen in the introduction. This relation is given by the following proposition, which is proved in \cite{Matmoun} following the lines of  \cite{Ga}.
\begin{prop}[see \cite{Matmoun}, Proposition 3.1]\label{GaMa}
Let $(M,g)$ be a pseudo-Riemannian manifold. 
Let  $(\hM,\hg)$ be the cone manifold over $(M,g)$.
There exists a smooth non-constant  function $\alpha: M\rightarrow \R$ such that for any vectorfields $X$, $Y$, $Z$ of $M$ we have:
$$(*)\qquad DDD \alpha (X,Y,Z)+ 2 (D\alpha  \otimes g)(X,Y,Z)+(D\alpha \otimes g)(Y,X,Z)+(D\alpha \otimes g)(Z,X,Y)=0,$$
if and only if there exists a non-trivial symmetric parallel $2$-tensorfield on $(\hM,\hg)$.

More precisely if $\alpha$ is a non-trivial solution of $(*)$ then the Hessian of the function $A:\hM\rightarrow \R$ defined by $A(r,m)=r^2\alpha(m)$ is parallel (ie $\hD\hD\hD A=0$). Conversely if $T$ is a symmetric parallel $2$-tensorfield on $\hM$ then $T(\partial_r,\partial_r)$ does not depend on $r$ and is a solution of $(*)$. Moreover $2T$ is the Hessian of the function $A_T$ defined by $A_T(r,m)=r^2T_{(r,m)}(\partial_r,\partial_r)$.
\end{prop}
We just quote the following Lemma, it is one of the steps of the proof of Proposition \ref{GaMa} and it will be useful further.
\begin{lem}[\cite{Matmoun}, Corollary 3.3]\label{retour}
 Let $T$ be a symmetric parallel $2$-tensorfield on $(\hM,\hg)$, and let $\alpha=T(\partial_r,\partial_r)$. Let $X$, $Y$, $Z$  be  vectors tangent to $M$ also seen as vectors perpendicular to $\partial_r$ in $\hM$. We have 
$$\begin{array}{rl}
   2T(\partial_r,X)=&rD\alpha(X)\\
   2T(X,Y)=&r^2(2g(X,Y)\alpha+DD\alpha(X,Y)),\\ 
2DT(X,Y,Z)=&  -D\alpha\otimes g (Y,X,Z) -  D\alpha\otimes g(Z,X,Y).
  \end{array}
$$
\end{lem}
\begin{example}\label{xamp}
 \rm{It follows from Fact \ref{fact}, that  the curvature of $\hg$ is given by 
$$\widehat R(X,Y)Z=R(X,Y)Z-g(Y,Z)X+g(X,Z)Y,$$ 
where $R$ and $\widehat R$ are the curvature tensors of $g$ and $\hg$.
Hence, a simply-connected pseudo-Riemannian manifold $(M,g)$  with constant curvature equal to $1$ is geodesically complete and has a cone which is flat and simply-connected. Hence the cone over $(M,g)$ admits any kind of parallel tensorfields. This is somehow a trivial example. We are going to look for non-trivial ones.
}
\end{example}
Contrarily to the Riemannian case, the existence of a parallel symmetric $2$-tensor on a pseudo-Rieman\-nian manifold does not imply that the manifold is decomposable (ie that it possess parallel non-degenerate distributions).
It is a consequence of the fact that the self-adjoint endomorphism associated to such a tensor and the metric can not always be simultaneously diagonalized. 
But the following Proposition shows that we can consider only three  cases.
\begin{prop}\label{above}
 If a pseudo-Riemannian manifold $(N,h)$ admits a non trivial symmetric parallel endomorphism  $\widetilde T$ then there exists on $(N,h)$ a  symmetric parallel endomorphism 
$\widetilde T'$ such that 
$\widetilde T'^2=\widetilde T'$, $\widetilde T'^2=0$ or $\widetilde T'^2=-Id$.
\end{prop}
{\bf Proof:} If $(N,h)$ is decomposable ie if there exists a non degenerate parallel distribution $V$ on $(N,h)$ then the projection on $V$ is a parallel endomorphism $\widetilde P$ satisfying $\widetilde P^2=\widetilde P$. If $\widetilde T$ is a 
 symmetric parallel endomorphism then it is also the case of its nilpotent part. If it is not trivial and if we take a proper power of it we obtain a non trivial symmetric parallel endomorphism  $\widetilde T'$ such that $\widetilde T'^2=0$.

At last, if $(N,h)$ is not decomposable and if $\widetilde T$ is a non trivial semi-simple symmetric parallel endomorphism then there exits $\lambda=a+ib  \in \C\setminus \R$ such that the minimal polynomial of $\widetilde T$ is $(X-\lambda)(X-\overline \lambda)$. But in this case $\widetilde T'=\frac{1}{b}\widetilde T-\frac{a}{b}\mathrm {Id}$ is parallel symmetric and satisfies $\widetilde T'^2=-\mathrm{Id}$. $\Box$
\hip
Example \ref{xamp} shows that any of these situations may occur on a cone over a complete pseudo-Riemannian manifold. However it is  proven in \cite{Matmoun} that on a cone over a \emph{compact} manifold there is only one type of symmetric $2$-tensorfield to investigate:
\begin{prop}[see \cite{Matmoun}, Proposition 3.4]\label{para->decomp}
 Let $(M,g)$ be a closed pseudo-Riemannian manifold. If the equation $(*)$ has a non-constant solution then $(\hM,\hg)$ is decomposable.
\end{prop}
\subsection{Basic properties.}
 Let $(M,g)$ be a pseudo-Riemannian manifold such that  its cone $(\hM,\hg)$ admits a non trivial symmetric parallel $2$-tensorfield $T$. Proposition \ref{above} tells us that we can assume that its associated parallel endomorphism $\widetilde T$ satisfies $\widetilde T^2=\widetilde T$, $\widetilde T^2=0$ or $\widetilde T^2=-Id$. 

We recall that if  $\lambda$ is an eigenvalue of $\widetilde T$ (here we can only have $\lambda=0\ \mrm{or}\ 1$) then the eigenspace $V_\lambda$ associated to $\lambda$ is a parallel distribution and therefore is integrable ie it defines a foliation. 


As in section \ref{section},
 we define on $\hM$  the functions $\alpha=T(\partial_r,\partial_r)$ and $A$ by $A(r,m)=r^2\alpha(m)$. We recall that $\alpha$ is actually a function on $M$. 

We start with two corollaries of Proposition \ref{GaMa}. The first actually implies Proposition \ref{para->decomp}.
\begin{cor}\label{critic}
The set of critical values of $\alpha$ is included in $\sigma(T)$, the spectrum of $\widetilde T$ (ie the set of its real eigenvalues). 
\end{cor}
{\bf Proof.} Let $m\in M$ be  a critical point of $\alpha$. 
 Lemma \ref{retour} and Proposition \ref{GaMa} implie that for any $r>0$,
$$\hD\hD A{}_{(r,m)}(\partial_r,.)=2\alpha(m)g(\partial_r,.).$$
 It means that $\partial_r(r,m)$ belongs to the eigenspace of  $\widetilde T$ associated to the eigenvalue $\alpha(m)$.
In our cases it implies that $\alpha(m)=0$ or $1$.
%
$\Box$
\begin{cor}\label{coco}
 If $\alpha$ is a solution of $(*)$ which is constant on an open subset $U$ of $M$ then $\alpha$ is constant on $M$.
\end{cor}
{\bf Proof.} As for any $k\in \R$, the function $\alpha+k$ is also a solution of $(*)$, we can assume that for any $m\in U$, $\alpha(m)=0$. As $D \alpha(m)=0$, it follows from 
Lemma \ref{retour}
that, for any $r>0$ and any $m\in U$, $\hD\hD A_{(r,m)} (\partial_r,.)$ vanishes on $TM=\partial_r^\perp$ and takes the value $2\alpha(m)$ on $\partial_r$. It means that
\begin{equation*}\hD\hD A_{(r,m)} (\partial_r,.)=2\alpha(m)g_{(r,m)}(\partial_r,.)=0.\label{decadix}\end{equation*}
Moreover as for any $m\in U$, we have $DD\alpha(m)=0$ then Proposition \ref{GaMa} and Lemma \ref{retour}
 give:
$$\hD\hD A_{(r,m)}(Y,Z)=2T_{(r,m)}(Y,Z)=2g(Y,Z)\alpha(m)=0.$$
Hence the Hessian of $A$ vanishes on $\R^*_+\times U$, as it is parallel, it vanishes everywhere. The gradient of $A$ is therefore parallel but it vanishes also on $\R^*_+\times U$ therefore $\alpha$ is constant.$\Box$
\hip

\hop
We consider the vectorfield $Y$ on $\hM$ defined by 
$$Y=\widetilde T (\partial_r).$$ 
We decompose now the vectorfield  $Y$ according to the splitting $T\hM=\R \partial_r \oplus TM$, we have 
$$Y=\alpha \partial_r + X,$$ 
 where $X$ is a vectorfield on $\hM$ tangent to $M$. We have:
$$\hg(X,X)=\left\{\begin{array}{ll}
\alpha-\alpha^2 &\quad \mrm{if}\quad \widetilde T^2=\widetilde T\\
-\alpha^2 &\quad \mrm{if}\quad \widetilde T^2= 0\\
-1-\alpha^2 &\quad \mrm{if} \quad \widetilde T^2= -Id
\end{array}\right.
$$
The following  proposition generalize  Corollary 4.1 of \cite{Leist} (which concerns the case $\widetilde T^2=\widetilde T$).
\begin{prop}\label{leist}
The vectorfield $2rX$ projects on a vectorfield on $M$ which is the gradient of $\alpha$ (with respect to the metric $g$). 
\end{prop}
{\bf Proof.} Let $(r,m)$ be a point of $\hM$. Let $Z$ be the lift of  vectorfield of $M$ perpendicular at $(r,m)$ to $X$. We known from proposition \ref{GaMa} that $\hD\hD A=2T$ therefore, using Lemma \ref{retour}, at $m$ we have:
$$D\alpha(Z)=\hD\hD A(\frac{1}{r}Z,\partial_r)=2\hg(\frac{1}{r}Z,Y)=0.$$
This means that $Z$ is perpendicular at $m$ to the gradient of $\alpha$ and therefore that $X$ projects on a well-defined direction field of $M$. As $g(2rX,2rX)$ does not depend on $r$ and, according to Corollary \ref{coco}, vanishes on a closed set with empty interior, the vectorfield $2rX$ does project on $M$.

To conclude we just have to show that $D\alpha(2rX)=g(2rX,2rX)$. Using again Lemma \ref{retour} we have:
$$D\alpha(2rX)=\hD\hD A(2X,\partial_r)=\hg(2X,Y)= 2\hg(2X,\alpha \partial_r + X)=4\hg(X,X)=g(2rX,2rX).$$
It means that the projection of  $2rX$ and the gradient of $\alpha$ coincide on a dense open set and therefore everywhere.
 $\Box$

\begin{cor}\label{gradient}
 The gradient of $A$ is the vectorfield $2rY$. It satisfies 
$$\hD_{rY}rY=\left\{\begin{array}{ll}
rY  &\quad \mrm{if}\quad \widetilde T^2=\widetilde T\\
0 &\quad \mrm{if}\quad \widetilde T^2= 0\\
-r\partial_r &\quad \mrm{if} \quad \widetilde T^2= -Id
\end{array}\right.
$$
and is therefore a pregeodesic vectorfield (ie up to reparameterization its integral curves are geodesics). The vectorfield $rX$ is pregeodesic for the metric $g$, more precisely we have 
$$D_{rX}rX=\left\{\begin{array}{ll}
(1-2\alpha)rX &\quad \mrm{if}\quad \widetilde T^2=\widetilde T\\
-2\alpha rX  &\quad \mrm{if}\quad \widetilde T^2= 0\\
-2\alpha rX  &\quad \mrm{if} \quad \widetilde T^2= -Id
\end{array}\right.
$$
\end{cor}
{\bf Proof.} We have $dA=2r\alpha dr+r^2 d\alpha$. Let $v=a\partial_r+h$ be a vector tangent to $\hM$ decomposed according to the splitting $T\hM=\R \partial_r \oplus TM $. We verify that $\hg(2rY,.)=dA$. Using proposition \ref{leist}, we have 
$$\begin{array}{ll}
dA(v)&=2r\alpha a+r^2d\alpha(h)\\
      &= 2r\alpha a+r^2 g(2rX,h)\\
      &= 2r \hg(\alpha \partial_r+X, v).
  \end{array} $$
The fact that the covariant derivative commutes with the musical isomorphisms ($\sharp$ and $\flat$) implies 
 $\hg(\hD_{\,.\,}2rY,\,.)$ is equal to the Hessian of $A$. Hence $\hg(\hD_{rY}2rY,\, .)=2T(rY,\, .)$
 and $\hD_{rY}rY=r \widetilde T^2 (\partial_r)$.

According to Fact \ref{fact} $D_{rX}rX$ is the projection on $TM$ of $\hD _{rX}rX$. Using the fact that $rX=rY-\alpha\partial_r$ and Fact \ref{fact} again, it is straightforward to compute $\hD_{rX}rX$ and therefore $D_{rX}rX$.
$\Box$ 
\begin{cor}\label{constant}
 The function $A$  is constant along the the leaves of the foliation spanned by $\ker \widetilde T$.
\end{cor}
{\bf Proof.}  It follows from the fact that for any $Z\in \ker (\widetilde T)$ we have $0=T(\partial_r, Z)=\hg(Y,Z)=dA(Z)$.$\Box$

\section{Description of the regular locus of $\alpha$}
\begin{defi}
 Let $\alpha$ be a  non trivial solution of equation $(*)$. We call \emph{regular locus} of $\alpha$ the open dense set of $M$ defined by $M\setminus\alpha^{-1}(\sigma(T))$, where $\sigma(T)$ denotes the spectrum of $\widetilde T$  (hence in our case $\sigma(T)\subset\{0,1\}$).
\end{defi}
We are able now to give a complete description of the regular locus of a solution $\alpha$. Everything starts from the following consequence of Section \ref{section}:
\begin{cor}\label{send}
 Let $U$ be the open dense subset of $M$ defined by $U=\{m\in M\, | \, g(rX,rX)\neq 0\}$. The  vectorfield   $\overline X$  defined on $U$ by $\overline X=-\frac{1}{\sqrt{|g(2rX,2rX)|}}2rX$ is geodesic (ie it satisfies $D_{\overline X}\overline X=0$) and its local flow preserves the foliation of $U$ by level sets of $\alpha$. Moreover, if $\gamma$ is an integral curve of $\overline X$, there exists a constant $c$ such that  we have 
$$\alpha(\gamma(s))=\left\{\begin{array}{ll}
\cos^2(s+c) &\quad \mrm{if}\quad \widetilde T^2=\widetilde T\quad  \mrm{and}\quad 0< \alpha < 1 \\
\cosh^2(s+c) &\quad \mrm{if}\quad \widetilde T^2=\widetilde T\quad  \mrm{and}\quad \alpha > 1 \\
-\sinh^2(s+c) &\quad \mrm{if}\quad \widetilde T^2=\widetilde T\quad  \mrm{and}\quad \alpha < 0 \\
e^{2t+c} &\quad \mrm{if}\quad \widetilde T^2= 0\quad  \mrm{and}\quad \alpha > 0\\
-e^{-2t+c} &\quad \mrm{if}\quad \widetilde T^2= 0\quad  \mrm{and}\quad \alpha < 0\\
\sinh(2s + c) &\quad \mrm{if} \quad \widetilde T^2= -Id
\end{array}\right.$$
\end{cor}
{\bf Proof:} 
The vectorfield $\overline X$ is pregeodesic and unitary, therefore it is geodesic. To see that it preserves the foliation of $U$ by level sets of $\alpha$, we compute $\overline X.\alpha$. According to Proposition \ref{leist} $2rX.\alpha=g(2rX,2rX)$, therefore we have
$$-\frac{1}{\sqrt{|g(2rX,2rX)|}}2rX.\alpha= -\frac{1}{\sqrt{|g(2rX,2rX)|}}g(2rX,2rX)= \left\{\begin{array}{ll}
\pm 2\sqrt{|\alpha-\alpha^2|}           &\quad \mrm{if}\quad \widetilde T^2=\widetilde T\\
2|\alpha|                             &\quad \mrm{if}\quad \widetilde T^2= 0\\
2\sqrt{1+\alpha^2}                   &\quad \mrm{if} \quad \widetilde T^2= -Id
\end{array}\right.
$$
The existence of such an equation implies that the local flow of  $\overline X$ preserves the foliation. Solving these ordinary differential equations, we obtain we expression of $\alpha$ along an integral line of $\overline X$.
$\Box$

Now, we are in order to state our first Theorem, its statement is quite long but a shorter version was given in the introduction.
\begin{thm}\label{lieuregulier}
Let $(M,g)$ be a pseudo-Riemannian manifold  such that $(\hM,\hg)$ has a parallel symmetric $2$ tensorfield $T$ such that $\widetilde T^2=\widetilde T$, $0$ or $-\mrm{Id}$. Let $\alpha$ be the function on $M$ associated to $T$ and $\sigma(T)$ be the spectrum $\widetilde T$.  Let $U$ be a connected component of $M\setminus \alpha^{-1}(\sigma(T))$.

 If $(M,g)$ is complete or if $\alpha$ has compact levels then there exist an open interval $I$, a pseudo-Riemannian manifold $(N,h)$ endowed with a parallel symmetric $2$ tensorfield $S$ (which can be trivial) such that $U=I\times N$ and
\begin{enumerate}
 \item\label{nil}
If $\widetilde T^2=0$, then $\widetilde S^2=0$, $\alpha(s,n)=e^{2s}$ and 
$$ 
g=-ds^2 + g_s=-ds^2 + (e^{2s}(h-S)+S).
$$
Each submanifold $(\{s\}\times N,g_s)$ is  endowed with the parallel symmetric $2$ tensorfield $S_s$ given by 
$$
S_s=e^{2s}S.
$$
Moreover if $(M,g)$ is non extendable (ie if it cannot be isometrically embedded in a manifold having the same dimension) then $I=\R$.
\item 
$\widetilde T^2=-\mrm{Id}$, then  $\widetilde S^2=-\mrm{Id}$, $\alpha(s,n)=e^{-\sinh^2(2s)}$ and 
$$ 
g=-ds^2+ g_s=-ds^2 +  h - \sinh(2s) S,
$$
Each  submanifold $(\{s\}\times N,g_s)$ is  endowed with the parallel symmetric $2$ tensorfield $S_s$ given by 
$$
S_s= S+\sinh(2s) h
$$
Moreover $U=M$ and if $(M,g)$ is non extendable then $I=\R$.
\item \label{f1}
If $\widetilde T^2=\widetilde T$ and $0<\alpha<1$,  then $\widetilde S^2=\widetilde S$,  $\alpha(s,n)=\cos^2(s)$ and 
$$ 
g=ds^2+ g_s=ds^2+ (\sin^2(s) (h-S) + \cos^2(s) S).
$$
Each submanifold $(\{s\}\times N,g_s)$ is  endowed with the parallel symmetric $2$ tensorfield $S_s$ given by 
$$
S_s=\cos^2 (s) S.
$$
Moreover if $(M,g)$ is non extendable then it is possible to choose $I=]0,\pi/2[$.
\item \label{f2}
If $\widetilde T^2=\widetilde T$ and $\alpha<0$ (or $\alpha>1$),  then $\widetilde S^2=\widetilde S$,  $\alpha(s,n)=\cosh^2{2s}$ (or $\alpha(s,n)=-\sinh^2(2s)$) and 
$$ 
g=-ds^2+ g_s=-ds^2 + (\sinh^2(s) (h-S) + \cosh^2(s) S).
$$
Each submanifold $(\{s\}\times N,g_s)$ is  endowed with the parallel symmetric $2$ tensorfield $S_s$ given by 
$$
S_s=\cosh^2 (s) S.
$$
Moreover if $(M,g)$ is non extendable then $I=\R_{>0}$.
\end{enumerate}
Reciprocally, if $(N,h)$ is a pseudo-Riemannian manifold endowed with a parallel symmetric $2$-tensor $S$ such that $\widetilde S^2=\widetilde S,\ 0$ or $-\mathrm{Id}$ then the manifolds $(I\times N, g)$ given by the above formulas are well defined and the cone over any of them admits a non-trivial parallel symmetric $2$ tensorfield $T$.
\end{thm}
{\bf Proof :} We begin the proof by the following Lemma.
\begin{lem}\label{un}
 Under the hypothesis of Theorem \ref{lieuregulier}, and keeping the notation therein, there exist 
  a manifold $N$ such that $U=\R\times N$. 
Moreover 
 each submanifold $\{x\}\times N$ correspond to a level of $\alpha|_{U}$ and the vectorfield $rX$ is everywhere tangent to the factor $\R$.
\end{lem}
{\bf Proof of Lemma \ref{un}:} Let $a$ and $b$ in $\R$ such that $a<b$ and $[a,b]\subset \alpha(U)$. 
Let $\varepsilon >0$ such that $[a-\varepsilon ,b+\varepsilon]\subset \alpha(U)$.
We denote by $K$ the subset of $U$ given by $K=\alpha^{-1}([a,b])$. We choose a bump function $l$ such that $l(K)=\{1\}$ and $l$ vanishes outside $\alpha^{-1}(]a-\varepsilon,b+\varepsilon[)$. As $g$ is  geodesically  complete or as $\alpha$ has compact levels, the vector field $l\overline X$ is complete.
 
Moreover Corollary \ref{send} says that the local flow  of $\overline {X}$ preserves the foliation of $U$ by level set of $\alpha$. It entails  that the flow of $l\overline X$ has the same property on $K$. We obtain that $\alpha^{-1}([a,b])$ is diffeomorphic to $J\times \alpha^{-1}(\{a\})$, where $J$ is a closed interval.

Thus $U$ is a fibre bundle over a $1$-dimensional manifold with fiber $N$ diffeomorphic to $\alpha^{-1}(\{a\})$. 
As $\alpha$ is constant along the fibers the base can not be compact, thus it is diffeomorphic to $\R$. The base being contractible,  the fibre bundle is trivial ie $U$ is diffeomorphic of $\R\times N$.
$\Box$

On the factor $\R$ given by Lemma \ref{un}, we choose a parameterization $t$ such that $\partial_t=rX$. 
The line is therefore identified to an open interval $J$
\begin{lem}\label{deux}
We denote by $\pi$ the projection on $J$, by $S_t$ and $g_t$ the restrictions of  $T$ and $g$  to $\{t\}\times N$, then
\begin{eqnarray}
\label{evolg} 
g'_t&=& -2\overline \alpha(t) g_t + 2  S_t\\
\label{evolS} 
S'_t&=& -2 \overline \alpha(t) S_t + 2 S^2_t,
\end{eqnarray}
where $\overline \alpha=\alpha\circ \pi$.
Moreover if we denote by $D^t$ the Levi-Civita connexion of $g_t$. We have 
$D^tS_t=0$.
\end{lem}
{\bf Proof of Lemma \ref{deux}:}
If we denote by $\mathcal L_{rX}$ the Lie derivative according to $rX$, we have
$$\mathcal L_{rX}g(u,v)=g(D_{u}rX,v)+g(D_vrX,v)=D^2\alpha(u,v).$$
By Lemma \ref{retour} we thus have
$$\mathcal L_{rX}g=-2\alpha g+2T,$$
which clearly implies $(\ref{evolg})$.

On the other hand, by Lemma \ref{retour} we have 
$$DT(rX,u,v)=0$$
for any $u$, $v$ perpendicular to $rX$. We deduce that if $u$ and $v$ are perpendicular to $rX$ then 
$$\begin{array}{rl}
\mathcal L_{rX}T(u,v)=& T(D_{u}rX,v)+T(D_vrX,u)\\
=& g(D_{u}rX,\widetilde T(v))+g(D_vrX,\widetilde T(u))\\
=& \frac{1}{2}D^2\alpha(u,\widetilde T(v))+\frac{1}{2}D^2\alpha(v,\widetilde T(u))\\
=& -\alpha g(u,\widetilde T(v))+T(u,\widetilde T(v))-\alpha g(v,\widetilde T(u))+T(v,\widetilde T(u))\\
=& -2\alpha g(u,\widetilde T(v)) + 2g(u, \widetilde T^2(v))
\end{array}
$$
Showing the second assertion. 

Let $U$, $v$ and $w$ three vectors tangent to $\{t\}\times N$ at a point $(t,n)$. Using Lemma \ref{retour} we can write
$$
0=DT(u,v,w)=D^tT_t(u,v,w)+T(II_t(u,v) \overline X, w)+T( II_t(u,w)\overline X, v),$$
where $II_t$ stands for the second fundamental form of $\{t\}\times N$.
But for any $n\in \N$, we have  $\widetilde T^n(\partial_r)\in \mrm{Span}(\partial_r,rX)$ therefore $T(\overline X, u)=T(\overline X, v)=0$. Consequently $D^tT_t=0$.
$\Box$

With our choice of parameter on $\R$ the metric $g$ has the following shape: $g=g(rX,rX)dt^2+g_t$. This not the desired one. Thus we need to reparameterize $\R$, we choose a parameter $s$ such that $\partial_s=\overline X$. The equations obtained at Lemma \ref{deux} turn into:\begin{eqnarray*}
\label{evolgs} g'_s&=& \frac{-1}{\sqrt{|g(rX,rX)}|}(-2\overline \alpha(s) g_s + 2  S_s)\\
\label{evolSs} S'_s&=&  \frac{-1}{\sqrt{|g(rX,rX)}|}(-2 \overline \alpha(s) S_s + 2 S^2_s).
\end{eqnarray*}
If $\widetilde T^2=0$, according to Corollary \ref{send}, we obtain the equations 
\begin{eqnarray*}
 g'_s&=& \epsilon 2 g_s - 2 \, e^{-\epsilon s}  S_s\\
 S'_s&=& \epsilon 2 S_s,
\end{eqnarray*}
where $\epsilon$ is the sign of $\alpha$.\\
If $\widetilde T^2=-\mrm{Id}$, according to Corollary \ref{send}, we obtain the equations 
\begin{eqnarray*}
 g'_s&=&  2 \tanh (2s) g_s - \frac{2}{\cosh(2s)} S_s\\
 S'_s&=&  2 \tanh (2s) S_s +  \frac{2}{\cosh(2s)} g_s.
\end{eqnarray*}
If $\widetilde T^2=\widetilde T$ and $0<\alpha<1$, according to Corollary \ref{send}, we obtain the equations 
\begin{eqnarray*}
g'_s&=& \frac{2}{\tan (s)}g_s - \frac{2}{\cos (s)\sin (s) } S_s\\
S'_s&=& -2\tan (s) S_s.
\end{eqnarray*}
If $\widetilde T^2=\widetilde T$ and $\alpha<0$, according to Corollary \ref{send}, we obtain the equations 
\begin{eqnarray*}
g'_s&=&\frac{2}{\tanh (s)}g_s -\frac{2}{\cosh(s)\sinh(s)} S_s\\
S'_s&=& 2\tanh(s) S_s.
\end{eqnarray*}
The case $\widetilde T^2=\widetilde T$ and $\alpha>1$ is similar to the later.
\hop
 There exists a unique solution to each of those systems of differential equations with given initial data.
Now there is no difficulty to check that those solutions are the one given in the statement. 
This is not a surprise but we remark that the endomorphisms $\widetilde S_s$ do not depend of $s$.
Thus we proved the first half of the theorem. 

For the reciprocal, the first thing to check is that the metrics $g_s$ given in the statement are never degenerate. As it is a pointwise property, it is linear algebra. 
We denote by $I_k$ the identity matrix at order $k$, and by $N$ the matrix $\begin{pmatrix}
                                                                             0&1\\0&0
                                                                            \end{pmatrix}$.
Let $\widetilde S$ and $\widetilde g_s$ such that $S=h(\widetilde S.,.)$ and $g_s=h(\widetilde g_s.,.)$. 
It is easy to find a frame of $T_pN$ such that the matrix of $\widetilde S$ is given by
$$\begin{array}{cl}
 \begin{pmatrix}
 I_k&0\\
0&0
\end{pmatrix} \quad &
\mathrm{when\ } \widetilde S^2=\widetilde S,\\
\begin{pmatrix}
N \\
&\ddots&\\
&&N\\
&&&0
\end{pmatrix} \quad &
\mathrm{when\ } \widetilde S^2=0,\\
 \begin{pmatrix}
 0&-I_{\frac{n}{2}}\\
I_{\frac{n}{2}}&0
\end{pmatrix} \quad &
\mathrm{when\ } \widetilde S^2=-\mathrm{Id}.
\end{array}
$$ The matrix in those frames of $\widetilde g_s$ are now easy to write down. They are clearly non degenerate.

If there exists a parallel tensor $T$ on $\hM$ inducing $S$ on $N$, then there exists a vectorfield $Y$ such that the  endomorphism $\widetilde T$ is given by~:
$$\left\{\begin{array}{ll}
 \widetilde T(\partial_r)=Y\\
 \widetilde T(Y)=0, -\partial_r\ \mrm{or}\ Y \quad \mrm{(according\ to\ the\ relation\ between\ } \widetilde S\ \mrm{and}\ \widetilde S^2\mrm{)}\\
 \widetilde T(Z)=\widetilde S(Z) \ \quad \mrm{for\ any\ vector}\ Z\ \mrm{tangent\ to} \{s\}\times N,
         \end{array}\right.
$$
As we already know the function $\alpha$, it not difficult to give $Y$. It has to be:
$$Y=\alpha(s) \partial_r - \frac{1}{r}\sqrt{\beta(s)}\partial_s,$$
 where 
$$\beta(s)=\left\{
\begin{array}{ll}
|\alpha-\alpha^2| &\quad \mrm{if}\quad \widetilde T^2=\widetilde T\\
\alpha^2 &\quad \mrm{if}\quad \widetilde T^2= 0\\
1+\alpha^2 &\quad \mrm{if} \quad \widetilde T^2= -\mrm{Id} 
\end{array},
\right.$$
 We have only to check that the tensorfield $T$ on $\hM$ given in the statement is parallel. The main difficulty is solved by the following Lemma.
\begin{lem}\label{trois}
 Let $(g_t,S_t)$ be path of metrics and symmetric $2$-tensors on a manifold $N$ satisfying (\ref{evolg}) and (\ref{evolS}) and such that $D^{t_0}S_{t_0}=0$. Then for any $t$ we have $D^tS_t=0$.
\end{lem}
{\bf Proof of Lemma \ref{trois}:} The first point is that $D^{t_0} S_t$ and $D^{t_0} g_t$ are solutions of a differential equation. Indeed we have:
\begin{eqnarray}
(D^{t_0} g_t)'&=&D^{t_0}(g'_t)=-2\alpha(t) D^{t_0}( g_t) + 2D^{t_0} (T_t)\\
(D^{t_0} T_t)'&=&D^{t_0}(T'_t)=- 2\alpha(t) D^{t_0}( T_t) +2D^{t_0}(T_t^2) 
\end{eqnarray}
As $\widetilde T^2=0, \, -\mrm{Id}\, \mrm{or}, \widetilde T$, this equation is in fact linear.
But for $t=t_0$ we have $D^{t_0}g_{t_0}=0$ and $D^{t_0}T_{t_0}=0$ therefore for any $t$ we have $D^{t_0}g_t=0$ and $D^{t_0}T_t=0$. It means that $D^{t_0}=D^t$ and $D^tT_t=0$.
$\Box$
\\
We leave to the reader the last verifications.$\Box$
\hop
The parts \ref{f1} and \ref{f2} of Theorem \ref{lieuregulier} have already appeared in \cite{Leist} under a slightly different form. To obtain the former version, we just have to apply De Rham-Wu theorem, locally or globally, to the triplet $(N,h,S)$. 
However, in order to write down the global version, we need an assumption of completeness on $N$ (see \cite{PR}), except in the Riemannian or in the Lorentzian case it is not enough to suppose the manifold compact. It is undoubtedly a very interesting question to know if there exist decomposable compact pseudo-Riemannian manifolds whose universal cover is not a product. 
\begin{cor}\label{alaleist}
 Let $(M,g)$ be a pseudo-Riemannian manifold such that $(\hM,\hg)$ has a parallel symmetric $2$ tensorfield $T$ such that
 $\widetilde T^2=\widetilde T$ (ie $(M,g)$ has a decomposable cone).
Let $U$ be a connected component of $M\setminus \alpha^{-1}(\{0,1\})$ and $\widetilde U$ denote its universal cover.
If the foliation defined by the kernel of the restriction of $T$ to $M$ is geodesically complete (for example if $g$ is geodesically complete or if $\alpha$ has compact levels and $\ker T$ is spacelike) then there exists two pseudo-Riemannian manifolds $(N_1,h_1)$ and $(N_2,h_2)$ such that 
  $$\begin{array}{rl}&\widetilde U_0=I\times N_1\times N_2\\ \mrm{and}\ & g=dts^2 + \cos^2(s) h_1 +\sin^2(s) h_2\\ \mrm{or\ }&  g=-ds^2 +\cosh^2(s) h_1 +\sinh^2(s)h_2.\end{array}$$
\end{cor}

The case $S=0$ of part \ref{nil} of Theorem \ref{lieuregulier} is also in \cite{Leist}. 
\section{Decomposable cones.}\label{plouf}
We suppose now that  $(M,g)$ is complete or that $\alpha$ is proper. The assumption of properness is stronger than the assumption that $\alpha$ has compact level, it is mainly  a way to say that the levels $\alpha^{-1}(0)$ or $\alpha^{-1}(1)$ are not empty.
It follows from Theorem \ref{lieuregulier} that the image of $\alpha$ is a union of sets chosen among the following one: $]-\infty,0[,\,\{0\},\,]0,1[,\,\{1\},\,]1,+\infty[$. 
\begin{defi}
We denote by $\ff$ the foliation of $\hM$ spanned by $\ker \widetilde T$ and by $\mathcal G$ the foliation defined by Im~$\widetilde T$.
\end{defi}
A parallel distribution being integrable those foliations do exist. Moreover their roles are symmetric as it is always possible to replace $\widetilde T$ by $\mrm{Id}-\widetilde T$. Those foliations will play an important role in what follows.

 To begin with, we suppose the function $\alpha$ bounded and the metric $g$ complete. Actually,  Theorem \ref{bounded} is almost proven in \cite{Leist} following Gallot's proof from \cite{Ga}. The missing point (the fact that the metric is Riemannian) is in \cite{Matmoun}. But for the convenience of the reader, we recall its proof (except for the Lemmas \ref{geodesic} and \ref{plat} that are stated without proof). We modified the presentation of some arguments in order to make clear how it is possible to adapt it to the case where $\alpha$ is suppose proper.
\begin{thm}\label{bounded}
 Let  $(M,g)$ be a complete pseudo-Riemannian manifold with decomposable cone.
We have $0\leq \alpha(m)\leq 1$ for  all $m\in M$ 
if and only if $(M,g)$ is finitely covered by a Riemannian round sphere.
\end{thm}
{\bf Proof.} We suppose $(M,g)$ complete. Let $Y=\widetilde T(\partial_r)$, $\Gamma$ be the geodesic of $(\hM,\hg)$ starting from the point $(r_0,m_0)$ in the direction $-rY(r_0,m_0)$ and $\gamma$ the geodesic of $(M,g)$ starting from $m_0$ in the direction $-rX(m_0)$. 
Let us remark that $\Gamma'(t)$ is never lightlike and that, by Corollary \ref{gradient}, $Y(\Gamma(t))$ is always proportional to $\Gamma'(t)$.

As in \cite{Leist} and \cite{Ga} we prove that $\Gamma$ contains a point where $Y$ vanishes. 
\begin{lem}[see \cite{Leist}, section 5]\label{geodesic}
 We denote by $\alpha_0$ the number $\alpha(m_0)$. We have $$\begin{array}{ll}
\gamma(t)=(r(t),\gamma(f(t))),\ \mrm{with}\\ 
r(t)=\sqrt{(\alpha_0\,t +r_0)^2+(\alpha_0-\alpha^2_0)r_0^2t^2}\\
f(t)=\frac{1}{\sqrt{\alpha_0-\alpha^2_0}}\arctan \left(\frac{\sqrt{(\alpha_0-\alpha^2_0)r_0t}}{\alpha_0t+r_0}\right)\quad  \mrm{if}\quad \alpha_0-\alpha^2_0 > 0 \\
f(t)=\frac{1}{\sqrt{\alpha^2_0-\alpha_0}} \mrm{argtanh} \left(\frac{\sqrt{\alpha^2_0-\alpha_0 r_0\, t}}{\alpha_0\, t+r_0}\right)\quad  \mrm{if}\quad \alpha_0-\alpha^2_0< 0 
\end{array} $$
\end{lem}
The geodesic $\Gamma$ is thus defined on $[0,1]$. We  deduce from Corollary \ref{gradient} that 
$$A\circ \Gamma'(t)=-2r_0\sqrt{\alpha_0}\sqrt{A\circ \Gamma(t)}.$$ 
Therefore $A\circ \Gamma(t)=(-r_0^2\alpha_0 t^2 +r_0^2\alpha_0)$ and $A\circ \Gamma (1)=0$. It implies that $\hg(Y(\Gamma(1)),Y(\Gamma(1)))=0$. As $Y(\Gamma(1))$ can not be lightlike, it has to be zero.

The point $\Gamma(1)$ is a minimum for $A$ therefore its Hessian positive, but (cf. Proposition \ref{GaMa}) $\hD\hD A=2T$ and $T$ is a projector. This means that the restriction of $\hg$ to Im~$\widetilde T$ is Riemannian.

But, we can replace $\widetilde T$ by $\mrm{Id}-\widetilde T$ and repeat this proof. We will obtain  that the restriction of  $\hg$ to $\ker \widetilde T$ is also Riemannian. 

It is well known (and it follows from fact \ref{fact}) that the curvature of $\hg$ is given by 
$$\widehat R(X,Y)Z=R(X,Y)Z-g(Y,Z)X+g(X,Z)Y,$$ 
where $R$ and $\widehat R$ are the curvature of $g$ and $\hg$. It implies that $(\hM,\hg)$ is flat if and only if $(M,g)$ has constant curvature equal to $1$. 

To prove that $\hg$ is flat, we use the following Lemma from \cite{Ga} and \cite{Leist}.
\begin{lem}[see \cite{Ga} Lemma 3.2 or  \cite{Leist} Lemma 6.3]\label{plat}
 If $Y(r,m)=0$ then the leaf of $\mathcal G$ (the foliation spanned by Im~$\widetilde T$) containing $(r,m)$ is flat.
\end{lem}
We have proven that $\mathcal G$ is flat. In order to prove  that $(\hM,\hg)$ is flat, we have to show that $\ff$ is also flat. It is done by repeating the proof with the tensorfield $\hg-T$ instead of $T$. $\Box$
\hop
As we said at the beginning of the section, we are also interested in replacing the assumption of geodesic completeness by the assumption that $\alpha$ is proper. 
To adapt the proof above  to  this case, we just have to prove the following proposition:
\begin{prop}\label{moi}
Let $m\in M$ such that $\alpha(m)<1$ and $\alpha(m)\neq 0$ and let $G_{(r,m)}$ be the leaf of $\mathcal G$ containing the point $(r,m)$.
If $\alpha$ is proper 
then there exists a point $p$  in $G_{(r,m)}$ such that $Y(p)=0$.
\end{prop}
{\bf Proof.} The proof starts the same way.
We are  looking for a critical point of the restriction of  $A$ to $G_{(r,m)}$. Classically we follow (backward) the integral curves of the gradient. We note that the gradient of the restriction of $A$ is also $2rY$ because $Y$ is tangent to $\mathcal G$. 
Let $(r,\gamma) :\ ]a,b[\ \rightarrow G{(r,m)}$ be the maximal integral curve of $-rY$ such that $(r(0),\gamma(0))=(r,m)$.
\begin{lem}\label{minore}
 The image of the  restriction of $\gamma$ to $[0,b[$ lies in a compact set of $\hM$.
\end{lem}
{\bf Proof.} 
The value of $\alpha$ is bounded along this curve, as $\alpha$ is proper, 
we just have to show that $r(]a,0])$ is contained in a compact subset of $]0,+\infty[$.

We suppose $\alpha(m)>0$ (respectively $\alpha(m)<0$): we first remark that $-rY.r=-r\alpha\leq 0$ (resp. $\geq 0$) 
This implies that $$\forall t\in [0,b[,  r(t)\leq r(0)=r\ \mrm{(resp.}\ r(t)\geq r(0)\mrm{)}.$$ We thus have a upper bound (resp. lower bound)  on $r(t)$.

If we apply Corollary \ref{constant} to the parallel tensorfield $\hg-T$, we obtain that the function $r^2(1-\alpha)$ is constant along the leaves of $\mathcal G$.

Hence, if $(r',m')$ is a point of $G_{(r,m)}$ 
 we have 
$r'^2(1-\alpha(m'))=r^2(1-\alpha(m))\neq 0$. 
Moreover $1-\alpha(m')\leq 1$ (resp. $1-\alpha(m')\geq 1$) therefore 
$r'\geq r_0\sqrt {1-\alpha(m)}\ \mrm{(resp}\ r'\leq r_0\sqrt {1-\alpha(m)}\mrm{)}$. As for all $t\in ]a,b[$, we have $(r(t),\gamma(t))\in G_{(r,m)}$, this gives a lower bound (resp. a upper bound) for $r(t)$.
$\Box$
\hop
It follows from lemma \ref{minore} that there exists a sequence $(t_n)_{n\in\N}$ of points of $[0,b[$ converging to $b$ and such that the sequence  $(\gamma(t_n))_{n\in\N}$ converges in $\hM$ to a point $(r_\infty, m_\infty)$. Let $O$ be  a foliated neighborhood for $\mathcal G$  of $(r_\infty, m_\infty)$. There are two possibilities: either $(r_\infty, m_\infty)$ belongs to $G_{(r,m)}$ or the points $\gamma(t_n)$ belong to an infinite number of connected components of $O\cap G_{(r,m)}$ (called plaques). The last case implies that the leaf $G_{(r,m)}$ accumulates around  $(r_\infty, m_\infty)$. As the vector $\partial_r$ is never tangent to $G_{(r,m)}$ this is incompatible with the following straightforward consequence of corollary \ref{constant}.
\begin{fact}
 Let $m\in M$, if the set $\R^*_+\times \{m\}\cap G_{(r,m)}$ contains more than  one point then  $\R^*_+\times \{m\}\subset  G_{(r,m)}$.
\end{fact}
Hence  $(r_\infty, m_\infty)\in G_{(r,m)}$ and is therefore a critical point of the restriction of $A$ to $G_{(r,m)}$. It means that the gradient vanishes at this point ie that $Y(r_\infty, m_\infty)=0$.
$\Box$

Replacing the $T$ by $\hg-T$ we get:
\begin{cor}
 Let $m\in M$ such that $\alpha(m)>0$ and $\alpha(m)\neq 1$ and let $F_{(r,m)}$ be the leaf of $\ff$ containing the point $(r,m)$.
If $\alpha$ is proper then there exists a point $p$  in $F_{(r,m)}$ such that $Y(p)=\partial_r$.
\end{cor}
We can therefore replace the hypothesis ``$g$ is complete'' by ``$\alpha$ is proper'' and repeat the proof of Theorem \ref{bounded}. However, if $\alpha$ is bounded and proper then $M$ is actually compact. Therefore, we can use Proposition \ref{para->decomp} to improve the statement. We get the following statement wich is actually the main result of \cite{Matmoun}:
\begin{cor}[see \cite{Matmoun} Theorem 1]\label{mainmatmoun}
 If $(M,g)$ is closed (compact without boundary) and if its cone admit a non trivial parallel symmetric $2$-tensorfield then $(M,g)$ is finitely covered by a Riemannian round sphere.
\end{cor}
{\bf Proof.}  Proposition \ref{para->decomp} says that as $M$ is compact then $(\hM,\hg)$ is decomposable. 
The manifold $M$ being closed, there exists $(m_+,m_-)\in M^2$  such that $\alpha(m_+)=\max_{m\in M} \alpha(m)$ and  $\alpha(m_-)=\min_{m\in M} \alpha(m)$, therefore $d\alpha(m_{\pm})=0$. According to corollary \ref{critic}, the only critical values of $\alpha$ are $0$ and $1$ therefore $\alpha (m_-)=0$ and $\alpha(m_+)=1$. 

Then Proposition \ref{moi} enable us to adapt the proof of Theorem \ref{bounded}.
$\Box$

There is a shorter way to prove Corollary \ref{mainmatmoun}. It consists in proving first that $(M,g)$ is Riemannian and therefore  complete (or apply Gallot-Tano theorem). It is what is done in \cite{Matmoun}. Anyway, our purpose was rather to use Proposition \ref{moi} than to give a proof of Corollary \ref{mainmatmoun}.
\hop
The following result extend Theorem \ref{bounded}, together they say that, for now on, there are no unexpected examples. 
\begin{thm}\label{cas2}
Let $(M,g)$ be a pseudo-Riemannian manifold with decomposable cone and $\alpha$ be the associated solution of equation $(*)$.
If $(M,g)$ is complete or if $\alpha$ is proper and if there exists $m\in M$ such that $0<\alpha(m)<1$ 
then $(M,g)$ has constant curvature equal to $1$.
\end{thm}
{\bf Proof.} 
We assume $(M,g)$  complete but using Proposition \ref{moi}, it is easy to adapt the proof to the case where $\alpha$ is  proper.

The first step consists in repeating the proof of Theorem \ref{bounded}. Doing so we obtain that  the metric $\hg$ is flat on $\R^*_+\times\alpha^{-1}(]0,1[)$ (it was also proven in \cite{Leist}).

Let $p\in \hM$ be a point such that $Y(p)=0$ (resp. $Y(p)=\partial_r)$. We are going to see that the curvature vanishes at $p$.
Corollary 
\ref{coco} tells us that
 $p$ belongs to the closure of $\R^*_+\times\alpha^{-1}(]0,1[)$ (and the curvature of $\hg$ therefore vanishes at $p$) unless   
$p$ is a local maximum (resp. minimum). 
But in that case the restriction of $\hg$ to Im~$\widetilde T$ (resp. $\ker \widetilde T$) is negative Riemannian. Hence $Y(p)$ (resp $Y(p)-\partial_r$) is never lightlike. Therefore if $\alpha(m)=0$ (resp $\alpha(m)=1$) then $Y(r,m)=0$ (resp. $Y(p)=\partial_r$). It means that for all $m\in M$ we have $\alpha(m) < 0$ (resp.$\alpha(m) > 1$) and this contradicts our hypothesis.

We suppose there exists a point $m\in M$ such that  $\alpha(m)< 0$ and we choose  $r>0$, then according to Lemma \ref{geodesic}
the geodesic starting from $(r,m)$ with initial speed $-rY$ (which is contained in a leaf of $\mathcal G$) is defined on $[0,1]$. Hence, reproducing the proof above, it contains a point $p$ such that $Y(p)=0$, then it follows from Lemma \ref{plat} 
that the curvature of $\hg$ vanishes along the leaf of $\mathcal G$ containing $(r,m)$. But as we just saw, the curvature of $\hg$  vanishes at $p$.
We have two  points $p$ and $(r,m)$ that lie in the same leaf of $\mathcal G$. As the distribution $\ker \widetilde T$ 
is  parallel and perpendicular to $\mathcal G$,  the curvature of the restriction of $\hg$ to $\ker \widetilde T$ is the same at $p$ and at $(r,m)$. It proves that $\hg$ is flat at $(r,m)$, therefore that $\hg$ is flat on $\alpha^{-1}(]-\infty,0])$.

To study the set $\alpha^{-1}([1,+\infty[)$, we consider the endomorphism $\mrm{Id}-\widetilde T$. The function associated to it is $1-\alpha$, hence this case is similar to the former.
We proved that $\hg$ is flat therefore that $g$ has constant curvature equal to $1$.$\Box$
\hip
The last case need more work, in particular we need a better understanding of the set $\alpha^{-1}(\{0,1\})$. It is also more interesting, as it provides examples with non constant curvature.
\begin{thm}\label{cas3}
 Let $(M,g)$ be a pseudo-Riemannian manifold with decomposable cone and $\alpha$ be the associated solution of equation $(*)$. We suppose $g$ is complete
or $\alpha$ is proper.
If for all $m\in M$, $\alpha(m)\leq 0$ or $\alpha(m)\geq 1$ then there exists a pseudo-Riemannian manifold $(N,h)$ such that
\begin{itemize}
 \item either  up to a $2$-cover, $M=\R\times N$ and  $g=-ds^2+\cosh ^2(s) h$.
\item or 
$\widetilde M$, the universal cover of $M$, is a warped product of the negative  $n$ dimensional hyperbolic space and the universal cover of $(N,h)$. \\
More precisely $\widetilde M$ is diffeomorphic to $\R^{n}\times \widetilde{N}$ and using polar coordinates on $\R^{n}$ the metric $g$ is given by $g=-ds^2+-\sinh^2(s)g_1+\cosh ^2(s) h$, where $g_1$ is the standard metric of the $(n-1)$-sphere.
Moreover 
 $M$ is a foliated bundle over $(N,h)$ with fiber $\H^{n}_-$, the $n$-dimensional negative Riemannian hyperbolic space, and with holonomy\footnote{this holonomy is the one that concerns foliations, it is not the pseudo-Riemannian one.} given by a morphism from the fundamental group of $N$  into the group of isometry of $\H^{n}$ fixing a point.
\end{itemize}
Moreover those manifolds are complete if and only if $N$ is complete and $\alpha$ is proper if and only if $N$ is compact.
\end{thm}
{\bf Proof.} If we replace $\widetilde T$ by $\mrm{Id}-\widetilde T$, we permute the cases  $\alpha(m)\geq 1$ and $\alpha(m)\leq 0$. Thus, without loss of generality,  we will now suppose that  for all $m\in M$, $\alpha(m)\leq 0$. We start the proof by the following lemma.
\begin{lem}\label{negative}
Let  $m$ be a point of $M$ and $\gamma$ be the geodesic such that $\gamma(0)=m$ and $\gamma'(0)=\overline X(m)$ then there exists $t\in\R$ such that $\alpha(\gamma(t))=0$.

The level set $\alpha^{-1}(0)$ is a $\dim(\ker(\widetilde T))-1$ dimensional submanifold. The restriction of $\hg$ to Im~$\widetilde T$ is negative definite.
\end{lem}
{\bf Proof.} 
The first point is given by Proposition \ref{moi} when $\alpha$ is proper and by Lemma \ref{geodesic} and the discussion that follows when $g$ is complete.

As $0$ is the maximum of $\alpha$, any element of  $\alpha^{-1}(0)$ is a critical point. As in Theorem \ref{bounded}, it entails that the restriction of $\hg$ to  Im~$\widetilde T$  is negative definite.

Furthermore, it follows from corollary \ref{constant} that if $\alpha(m)=0$ then $A$ vanishes on any leaf of $\ff$ (the foliation spanned by $\ker \widetilde T$) containing a point $(r,m)$. 
The Hessian of $A$ is twice the restriction of $\hg$ to Im~$\widetilde T$, hence 
 the singular points of the restriction of $A$ to the leaves of $\mathcal G$ are isolated. It means that $A^{-1}(0)$ is given by a reunion of isolated leaves of $\ff$. Moreover the vectorfield $\partial_r$ is everywhere tangent to those leaves. The set  $\alpha^{-1}(0)$ being the projection of $A^{-1}(0)$ is therefore a $\dim(\ker(\widetilde T))-1$ dimensional smooth submanifold. 
$\Box$
\begin{lem}\label{plan}
 The projection of the distribution Im~$\widetilde T$ to $M$ is a smooth  integrable totally geodesic timelike distribution. We denote it by $V$ and by $\mathcal G'$ the foliation it defines. 
If $n_1= \dim \mrm{Im}\ \widetilde T>1$, then any  leaf of $\mathcal G'$ is isometric to the negative definite hyperbolic space $\H^{n_1}_-$
\end{lem}
{\bf Proof.} As $\partial_r$ is geodesic and Im~$\widetilde T$ is parallel, we know that the projection of Im~$\widetilde T(r,m)$ on $T_mM$ does not depend on $r$. Moreover Im~$\widetilde T$ never contains $\partial_r$ therefore  the distribution   $V$ is  smooth and integrable. 
The restriction of $\hg$ to Im~$\widetilde T$ is negative definite, therefore the restriction of $g$ to $V$ is also negative definite.

Let $Z,Z'$ be two vectorfields  tangent to $V$. Their lift to $\hM$, still denoted by $Z,Z'$ lie in  Im~$\widetilde T\oplus \R\partial_r$. This last distribution is clearly  totally geodesic. Moreover, from fact \ref{fact}, we know  that $\hD_{Z}Z'={D_ZZ'}-rg(Z,Z')\partial_r$. It means that $D_ZZ'$ is tangent to $V$ and therefore that $V$ is totally geodesic.

Let $u$ and $v$ be two vector of $T_mM$, from Lemma \ref{retour} we have  
\begin{equation}2T(u,v)=\hD\hD A(u,v)=
r^2 (DD\alpha(u,v)+2 g(u,v)\alpha(m)).\label{hess}\end{equation}
If $m$ is a critical point of $\alpha$, we have $\alpha(m)=0$ and Im~$\widetilde T$ perpendicular $\partial_r$.
Hence the Hessian of $\alpha$ at a critical point is given by the restriction of $g$ to $V$. Therefore the restriction of $\alpha$ to any leaf $G'$ of $\mathcal G'$ is a Morse function and the critical points of $\alpha|_{G'}$ are isolated. 
If $n_1>1$ the set of regular points of the restriction of  $\alpha$ to $G'$ is connected. 

Moreover as the vectorfield $2rX$ is tangent to $\mathcal G'$ it is also the gradient of $\alpha|_{G'}$. Hence (using geodesic completeness or the properness of $\alpha$) there exists a critical point of $\alpha$ in the closure of $G'$. But $\alpha^{-1}(0)$ is transverse to $\mathcal G'$, therefore $\alpha$ vanishes on $G'$.

We use now the vector field $\overline X=-rX/\sqrt{-g(rX,rX)}$, which is defined on $M\setminus \alpha^{-1}(0)$. 
 Let $c_1<c_2<0$ 
and 
$\varepsilon>0$ such that $c_2+\varepsilon <0$.
Let $l$ be a function vanishing outside $\alpha^{-1}(]c_1-\varepsilon,c_2+\varepsilon[)$ and being constant equal to one on $\alpha^{-1}(]c_1,c_2[)$. The vectorfield $l\,\overline X$ is complete in both cases. 
According to Corollary \ref{send}, its  flow restricted to $\alpha^{-1}(]c_1,c_2[)$ sends level sets of $\alpha$ on level sets of $\alpha$. It follows that any two regular level sets of  $\alpha|_{G'}$ are diffeomorphic. 

If $\alpha|_{G'}^{-1}(c_1)$ is not connected, we can saturate its connected component by the gradient line in order to obtain a partition of the set of regular points of $G'$. This set being connected, we have a contradiction, therefore the level sets of $\alpha|_{G'}$ are connected. 

Hence, as we have a Morse function any level set of the restriction of $\alpha$ to a leaf of $\mathcal G'$ is diffeomorphic to the sphere $S^{n_1-1}$. Consequently the leaves of $\mathcal G'$ are diffeomorphic to $\R^{n_1}$.

We have shown that the function $\alpha$ vanishes on any leaf of $\mathcal G'$ therefore the function $A$ vanishes on any leaf of $\mathcal G$. It follows from Lemma \ref{plat} that the curvature of $\hg$ vanishes along Im~$\widetilde T$, and by Fact \ref{fact} that the curvature of $g$ is constant equal to one on $V$. 

Thus we know that a leaf of $\mathcal G'$ is diffeomorphic to $\R^{n_1}$, is negative definite and has constant curvature equal to $1$, as moreover the flow of $\overline X$ is future complete (even if we suppose $\alpha$ proper) we can say that any leaf of $\mathcal G'$ is isometric to $\H^{n_1}_-$.
$\Box$
\hop
We denote by $W$ the distribution of $TM$ which is the orthogonal complement of $V$. This distribution is integrable, its  leaves being the intersection of the level sets of $\alpha$ (which are all smooth) with the projection on $M$ of the leaves of $\ff$. We denote by $\ff'$ the foliation spanned by $W$. In particular $\alpha^{-1}(0)$ is a reunion of leaves of $\mathcal F'$.
\hop
We assume now  that $n_1=\dim \mrm{Im}~\widetilde T=1$. Taking eventually a $2$-cover, the direction field given by the projection of Im~$\widetilde T$ is well defined and \emph{oriented}. Moreover the vector field $2rX$ belongs to this field. It means that it is possible to extend to $M$ the vectorfield  $\overline X=rX/\sqrt{-g(rX,rX)}$. 
The function $\alpha$ being unbounded by Corollary \ref{send} this vectorfield 
has to be complete.
Always by Corollary \ref{send} it sends levels of $\alpha$ on levels of $\alpha$. Moreover the leaves of $\mathcal G'$ cannot be compact, therefore there exists a manifold $N$ such that $M$ is diffeomorphic to $\R\times N$.
At last, Corollary \ref{alaleist} (see also theorem 7.1 from \cite {Leist}) shows that $g$ is given on $\R^*\times N$ by 
$-ds^2+\cosh ^2(s) g_2$. The first assertion of the Theorem follows by continuity.
\hop
We assume now that  $n_1=\dim \mrm{Im}~\widetilde T>1$. 
\begin{fact}\label{fact2} 
The submanifold $\alpha^{-1}(0)$ is connected.  The universal cover of $M$ is diffeomorphic to $\R^{n_1}\times \widetilde N$, where $\widetilde N$ is the universal cover of $\alpha^{-1}(0)$.
\end{fact}
{\bf Proof.}
The foliation $\mathcal G'$ is totally geodesic and complete (it is isometric to $\mathbb H^{n_1}_-$), $M$ is simply connected, therefore by Theorem 2 of \cite{PR}, $\widetilde M$, the universal cover $M$ is diffeomorphic to the product of the universal cover of a leaf of $\mathcal G'$ by the  universal cover of a leaf of $\mathcal F'$. The function $\alpha$ vanishes only once on each leaf of $\mathcal G'$. Consequently $\alpha^{-1}(0)$ is connected and is equal to exactly  one leaf of $\mathcal F'$. $\Box$
\hip
We will denote by $ \widetilde{\alpha}$, $\widetilde{\ff'}$ and $\widetilde{\mathcal G'}$ the lifts to $\widetilde M$ of $\alpha$, $\ff'$ and $\mathcal G'$. 

We note that, as $n_1>1$,  $\widetilde M\setminus \widetilde {\alpha}^{-1}(0)$ is connected. Hence fact \ref{fact2} implies that $\widetilde M\setminus \widetilde {\alpha}^{-1}(0)$ is diffeomorphic to $\R_+^* \times S^{n_1-1}\times \widetilde N$, the factor $\R$ corresponding to the direction of the vectorfield $\overline X$. Furthermore by Theorem \ref{lieuregulier} 
the metric can be written as $-ds^2 - \sinh^2(s)g_1+\cosh^2(s) g_2$ (the coordinate $s$ is a priori defined up to a constant) 
where $g_1$ is a metric on $S^{n_1-1}$ and $g_2$ a metric on $\widetilde N$. We can see that $-ds^2+-\sinh^2(s)g_1$ gives the metric on $\mathcal G'$. Thus $ds^2 + \sinh^2(s)g_1$ corresponds to the metric of $\mathbb H^{n_1}$ in polar coordinate and $g_1$ is the canonical metric of the sphere.

The manifold  $M$ is the quotient of $\R^{n_1}\times \widetilde N$ endowed with this warped metric by the  action of   its fundamental group $\Lambda$. This action preserves in particular the metric, the function $\widetilde{\alpha}$ and the foliations $\widetilde{\mathcal G'}$ and $\widetilde{\mathcal F'}$. It means that $\Lambda=\Lambda_1\times \Lambda_2$, where $\Lambda_1\subset \mathrm{Isom}(\H^{n_1}_-)$ and $\Lambda_2\subset \mathrm{Diff}(N)$.
It is proven in lemma \ref{plan} that the leaves of $\mathcal G'$ are always diffeomorphic to $\R^{n_1}$. Thus an element  of $\Lambda$ fixing a point of $N$ has to be trivial. It means that there exists a morphism $\rho: \Lambda_2\rightarrow \mathrm{Isom}(\H^{n_1}_-)$ such that $\Lambda_1=\rho(\Lambda_2)$. Moreover the action $\Lambda$ on $M$ preserves $\widetilde {\alpha}^{-1}(0)$, therefore $\Lambda_1$ fixes a point, $\Lambda_2\subset\mathrm{Isom}(\widetilde N,g_2)$ and  $\Lambda_2$ acts properly discontinuously on $\widetilde N$ (and $\widetilde N/\Lambda_2=N$).  In that case, the action of $\rho(\Lambda_2)\times \Lambda_2$ is always properly discontinuous, isometric and it preserves the foliations.

To investigate geodesic completeness we will use the following proposition, where $f$  denotes the warping function.
\begin{prop}[see \cite{Oneill} Proposition 7.38 p.208]\label{Oneill}
A curve $\gamma=(\gamma_1,\gamma_2)$ in $M=\H^{n_1}_-\times_f N$ is a geodesic if and only if
\begin{eqnarray}
 \gamma_1''&=&h(\gamma'_2,\gamma'_2)f\circ\gamma_1\ \mrm{grad} f\\
\gamma_2'' &=& \frac{-2}{f\circ \gamma_1}(f\circ \gamma_1)'\gamma_2'
\end{eqnarray}
\end{prop}
Our function $f$ has a critical point $O$ (it verifies $f(O)=1$).
We deduce from  Proposition \ref{Oneill} that $\{O\}\times N$ is totally geodesic. 
Hence, if $(M,g)$ is geodesically then $(N,h)$ is also complete.

Reciprocally, we suppose $(N,h)$ complete. 
Let $\gamma(t)=(\gamma_1(t),\gamma_2(t))$ be a geodesic of $(M,g)$ and 
let $\Gamma_2$ be the locus of the geodesic of $(N,h)$ with initial speed $\gamma_2''(0)$. Let $u$ be a geodesic parameterization of $\Gamma_2$ (for example a arclength parameterization when the geodesic is not ligthlike).
It follows from Proposition \ref{Oneill}  that $\mathbb H^{n_1}_-\times \Gamma_2 $ is a totally geodesic submanifold of $(M,g)$.  There are three cases to consider according to the type of $\gamma_2'(t_0)$.

If $h(\gamma_2',\gamma'_2)=0$. In this case $\gamma_1$ is a geodesic of $\H^{n_1}_-$. We write $\gamma_2(t)$ as $\gamma_2(u(t))$. We have 
$$u''(t)= \frac{-2(f\circ \gamma_1)'(t)}{f\circ \gamma_1(t)}u'(t)$$
therefore $u'(t)=C(f\circ \gamma_1)^{-2}(t)$. Moreover $f\circ \gamma_1\rightarrow \infty$ implies $t\rightarrow\infty$. Hence the geodesic is complete.

If $h(\gamma_2',\gamma'_2)<0$ the restriction of $g$ to $\mathbb H^{n_1}_-\times \Gamma_2 $ is negative Riemannian and therefore complete (see \cite{Oneill} Lemma 7.40 p.209).

If $h(\gamma_2',\gamma'_2)>0$ the restriction of $g$ to $\mathbb H^{n_1}_-\times \Gamma_2 $ has constant curvature equal to $1$ (see for example Corollary 2.3 of \cite{Leist}). Moreover its signature is $(1,n_1)$ and it contains a codimension $1$ foliation by hyperbolic spaces and a complete geodesic (the one above $O$) perpendicular to this foliation. It means that the universal cover of  $\mathbb H^{n_1}_-\times \Gamma_2 $ is the universal cover of the (negative) anti de Sitter space and is therefore complete.

As any geodesic in $(M,g)$ is contained in such a submanifold, we have proven that $(M,g)$ is complete.
$\Box$
\section{Cones admitting a parallel nilpotent symmetric endomorphism field}\label{sect-nilpo}
To start the construction of the examples announced in the introduction, we take a manifold which is trivially an example: a complete, simply connected pseudo-Riemannian manifold with constant curvature equal to $1$. Its cone being flat and simply connected it admits any kind of parallel tensor.

We denote by $\R^{p+1,q}$ the space $\R^{p+q+1}$ equipped with  the standard pseudo-Euclidean metric of  signature  $(p+1, q)$. We choose coordinate such that the metric writes $2x_1x_2 + \sum x_i^2-\sum x_j^2$.  We consider the  pseudo-sphere $S^{p,q}=\{x\in \R^{p+1,q}\,|\, \langle x,x\rangle=1\}$.
The cone over $S^{p,q}$ is the open subspace $\R^{p+1,q}$ given by $\{x\in \R^{p+1,q}\,|\, \langle x,x\rangle>0\}$, the vectorfield $\partial_r$ being $\frac{1}{\sqrt{\langle x,x\rangle}}\,x$.

Obviously, the function on  $\R^{p+1,q}$ having a non trivial parallel Hessian are the polynomial function of degree $2$.
We choose $A(x)=x_1^2$. Its Hessian is defined by the $2$-step nilpotent endomorphism $2\widetilde T$ where 
$\widetilde T(u)=\hg(\partial_2,u)\partial_2$. The function $\alpha$ is just the restriction of $A$ to $S^{p,q}$ and 
$\alpha^{-1}(0)$ is clearly a codimension $1$ totally geodesic lightlike submanifold. Actually, the level set of $\alpha$ define a smooth codimension $1$ foliation ($\alpha$ is not a submersion but it has a square root which is a submersion).

Each connected component of  $S^{p,q}\setminus \alpha^{-1}(0)$ is isometric to $\R\times \R^{p+q-1}$ endowed with the metric $g_0=-dt^2+e^{-2t}h_0$ with $h_0$ flat.

It is not difficult to see that the gradient of $\alpha$ is given by $\mrm{grad}\ \alpha= 2rX=2x_1\partial_2-2x_1^2\partial_r$. Moreover, according to Theorem \ref{lieuregulier} we can choose the isometry in such a way that it sends the vectorfield $\partial _t$ to $\overline X=\frac{1}{\sqrt{\langle 2rX,2rX\rangle}}2rX$. 

It is not difficult to see that $\lim_{t\rightarrow +\infty}x_2(t)=+\infty$. Let $h_1$ be a perturbation with compact support of the metric $h_0$ on $\R^{p+q-1}$. We endow $S^{p,q}\setminus \alpha^{-1}(0)$ with the metric  $g_1=-dt^2+e^{-2t}h_1$. There exists an open neighborhood $U$ of $\alpha^{-1}(0)$ such that the metrics $g_0$ and $g_1$ coincide on $U\setminus \alpha^{-1}(0)$. Consequently the metric $g_1$ extend to a smooth metric $g$ on $S^{p,q}$. 

According to Theorem \ref{lieuregulier}, the restriction of $\alpha$ is a solution of $(*)$ on the manifold $S^{p,q}\setminus \alpha^{-1}(0)$ endowed with the  metric $g_1$. As this set is dense in $S^{p,q}$ and as $g_1$ can be extended to $S^{p,q}$, the function $\alpha$ is a solution on  of $(*)$ on the manifold $(S^{p,q},g)$.

We can construct this way a  lot of non flat manifolds whose cone admits a parallel nilpotent symmetric $2$ tensorfield. 

Those examples show clearly that the communication between $\alpha^{-1}(0)$ and the $M\setminus \alpha^{-1}(0)$ is not as simple as in the decomposable case. In particular, the gradient line, even if they are complete, never reach  critical points.

Let us look at the lack of those examples:
\begin{itemize}
\item We did not prove that $g_1$ can be chosen to be complete. The main reason is that they are clearly not extendable therefore, in some sense, complete enough. 
\item we have chosen a situation where the rank of $\widetilde T$ is $1$. It allowed us to perturb $h_0$ without thinking to $T$.  If $\widetilde S$ is a parallel nilpotent symmetric endomorphism on a flat manifold $(N,h)$, it is possible to perturb $h$ while keeping the endomorphism $\widetilde S$ parallel. Hence, it is possible to provide examples with more sophisticated $\alpha$ and $T$ but it is more technical.
\item Those examples are still flat near $\alpha^{-1}(0)$. This is perhaps the main problem. We did not start the discussion about all the admissible metric $h_1$. For example it seems reasonable to think that we can replace the metric $h_1$ by a metric asymptotically flat. In fact, we prefer ask if there exist non flat real analytic examples. Such an example would make pointless the discussion about the behavior at infinity of $h_1$.
\end{itemize}
\section{Application to projective geometry}
We define the \emph{degree of mobility} of a pseudo-Riemannian metric $g$ as the dimension of the space of metrics projectively (or geodesically) equivalent to $g$ ie the set of metrics having the same \emph{unparameterized} geodesics as $g$. This number is well defined see \cite{KM} for details. It is always positive as the connection is invariant when $g$ is multiplied by a constant. We will say that two metric are affinely equivalent if they have the same \emph{parameterized} geodesics ie their  Levi-Civita connections coincide.

Using results of \cite{KM},  Matveev and the author proved in \cite{Matmoun} the following result:
\begin{thm}
Let $g$ be a  pseudo-Riemannian metric on an $(n>1)-$dimensional  closed connected manifold. 
     Then,  if the metric $\bar g$ on $M$  is geodesically equivalent to $g$, but not affinely equivalent to $g$, then the degree of mobility  of $g$ is precisely $2$ or there exists $c\neq 0$ such that  $c\,g$ is Riemannian and has constant  curvature equal to one.
\end{thm}
We will see that this result is no more true when $g$ is complete and not compact.
The following result clarifies the link between the  Obata equation and projective geometry. It does not pretend to be new (see \cite{KM} for example). We tried to make clear how a metric is obtained from a parallel tensor $T$ on $\hM$.
\begin{prop}\label{projo}
 If $(M,g)$ is a manifold such that its cone admits a parallel symmetric $2$-tensorfield $T$ such that $T(\partial_r,\partial_r)$ is bounded above or below then there exists a metric $g'$ on $M$ which is projectively equivalent to $g$ but not geodesically equivalent to $g$.
\end{prop}
{\bf Proof:} Let $T$ be a symmetric $2$-tensorfield on $\hM$. As for any $(a,b)\in \R^2$ the tensor-field $aT+b\hg$ is parallel, we can assume that it is non degenerate and $\alpha=T(\partial_r,\partial_r)$ is positive. Thus it is a pseudo-Riemannian metric. Its Levi-Civita connection is the only torsion free connection $\nabla$ such that $\nabla T=0$. As $\hD T=0$ this connection is $\hD$.

According to fact \ref{fact}, we have 
$$\mathcal L_{r\partial_r} T(u,v)=T(D_u(r\partial_r),v)+ T(D_v(r\partial_r),u)=2T(u,v).$$
therefore, if we denote by  $T'$ the retriction of $T$ to the $T$-orthogonal of $\partial_r$, we have:
$$T=\alpha dr^2+  r^2 \frac{\alpha}{\alpha} T'$$
Stating $\rho=\sqrt{\alpha}r$ and $g'= \frac{1}{\alpha} T'$, we obtain
$$T=d\rho^2+\rho^2 g'.$$
It means that $T$ can also be seen as a cone metric and that $g'$ is obtained by taking the restriction of $T$ on the submanifold $M_\alpha$ of $\hM$ defined by  $M_\alpha=\{(\alpha(m),m)\in \hM \,|\, m\in M\}$.
\begin{lem}\label{samegeod} 
The projection of a geodesic of $(\hM,\hD)$ on $M$ is a pregeodesic of $(M,g)$ and of $(M,g')$.
\end{lem}
{\bf Proof of Lemma \ref{samegeod}:}
Let $\Gamma(t)=(r(t),\gamma(t))$ be a geodesic of $(\hM,\hD)$. In follows from Fact \ref{fact} that $r(t)$ and $\gamma(t)$ satisfy:
\begin{eqnarray*}
 0&=& \ddot{r}(t)-r(t)g(\dot\gamma(t),\dot\gamma(t))\\
0&=& 2\dot r(t)\dot \gamma (t) +r(t) D_{\dot\gamma}\dot\gamma(t)
\end{eqnarray*}
Hence  $ D_{\dot\gamma}\dot\gamma(t)$ is proportionnal to $\dot \gamma (t)$, it means that $\gamma(t)$ is a pregeodesic for $(M,g)$. Using parameter $\rho$ instead of $r$, the same proof can be done for $(M,g')$. $\Box$
\hop
We proved that the metrics are projectively equivalent. To see that they are not affinely equivalent, 
it is possible to compute explicitely the reparameterization of the geodesics. It can be done with the help of section 5 from \cite {Leist}, where the geodesic lift of a geodesic of $(M,g)$ to $\hM$ is explicitely computed (see also \cite{KM} section 2.4). There always exist geodesics that are not affinely reparameterized. 
$\Box$
\hop
\begin{example}\label{ex}
{\rm  
\begin{itemize}
 \item  Let $(M,g)=\H^{n}_-\times_f N$ be one of the  manifolds obtained in Theorem \ref{cas3}. Its cone has a $n$ dimenional negative Riemannian flat distribution. Let $\widetilde T_1$ be the projection on the tangent space of this distribution.  We have $T_1(\partial_r,\partial_r)= g(\widetilde T_1(\partial_r),\partial_r)\leq 0$ so we consider $T_1'=g-\frac{1}{2}T_1$. It gives us a metric $g'_1$ which is  projectively equivalent metric to $g$.

Moreover, if $n>1$, we can split the flat distribution in two. It gives two parallel distributions. Choosing one of them, we obtain this way an other endomorphism $\widetilde T_2$ and an other geodesically equivalent metric $g'_2$. Hence the degree of mobility of $g$ is greater than $2$.
\item We can also use the example of section \ref{sect-nilpo}, to obtain a certainly less known family of projectively equivalent metrics with non constant curvature. But in this case, the degree of mobility is probably equal to $2$.
\end{itemize}
}
\end{example}
We have proven:
\begin{cor}
 There exists  pseudo-Riemannian manifolds $g$ and $g'$ with non constant curvature such that 
\begin{enumerate}
\item $g$ is geodesically complete,
\item $g$ and $g'$ are projectively equivalent but not affinely equivalent,
\item the degree of mobility  of $g$ and $g'$ is greater than $2$.
\end{enumerate}
\end{cor}
It seems natural to  ask now if the metrics given in the example \ref{ex} are the only complete  pseudo-Riemannian manifolds with non constant curvature and a degree of mobility greater than $2$ that admit projectively but not affinely equivalent metrics. Indeed Kiosak and Matveev have shown in \cite{KM} that if $(M,g)$ is a  pseudo-Riemannian metric with a degree of mobility  is greater than $2$ then any metric $g'$, which is projectively equivalent but non affinely equivalent to $g$,   is associated to a non-trivial solution of the Obata equation. They also have shown that the metric obtained this way are never complete.
It follows from this result and our study that any other example should have a cone with a nilpotent parallel endomorphism.

\hip
\begin{tabular}{ll}
 Address: & Universit\'e Bordeaux 1, Institut de Math\'ematiques de Bordeaux,\\
 &351, cours de la lib\'eration, F-33405 Talence, France\\
E-mail:&{\tt pierre.mounoud@math.u-bordeaux1.fr}
\end{tabular}
\end{document}